\documentclass[12pt]{amsart}
\usepackage{epsfig}
\usepackage{amsmath}
\usepackage{amssymb}
\usepackage{amscd}
\usepackage{color}

\newtheorem{thm}{Theorem}[section]

\newtheorem{pro}[thm]{Proposition}
\theoremstyle{definition}

\newtheorem{remark}{Remark}

\numberwithin{equation}{section}

\begin{document}

\title[diffeomorphism with global dominated splitting can not be minimal]
{A diffeomorphism with global dominated splitting can not be
minimal}
\author[Pengfei Zhang]{Pengfei Zhang}

\address{Department of Mathematics, University of Science and Technology of China,
 Hefei, Anhui 230026, P. R. China\\
\newline and
CEMA, Central University of Finance and Economics, Beijing 100081,
China }



\email{pfzh311@gmail.com}

\begin{abstract}
Let $M$ be a closed manifold and $f$ be a diffeomorphism on $M$. We
show that if $f$ has a nontrivial dominated splitting $TM=E\oplus
F$, then $f$ can not be minimal. The proof mainly use Ma\~n\'e's
argument and Liao's selecting lemma.
\end{abstract}

\maketitle

\section{Introduction}

In \cite{H} Herman constructed a (family of) $C^\infty$
diffeomorphism on a compact manifold that is minimal and has
positive topological entropy simultaneously. So positive entropy is
insufficient to guarantee the nonminimality. This draws forth the
problem to find some nature structure of the system that
incompatible with the minimality. In \cite{M} Ma\~n\'e gave an
argument to locate some nonrecurrent point if the map admits some
invariant expanding foliation (also see \cite{BRHRHTU}). In
particular this argument shows that a partially hyperbolic
diffeomorphism always has some nonrecurrent point and hence can not
be minimal. In this note we show that a global dominated splitting
is sufficient to exclude the minimality of the system.

Let $M$ be a closed Riemannian manifold and $f:M\to M$ be a $C^1$
diffeomorphism on $M$. The map $f$ is said to have a global
dominated splitting on $M$ if there exist an invariant splitting
$TM=E\oplus F$, two numbers $\lambda\in(0,1)$ and $C\ge1$ such that
\begin{equation}\label{domi}
\|Df^n|_{E(x)}\|\cdot\|Df^{-n}|_{F(f^nx)}\|<C\lambda^n\text{ for all
}n\ge1,x\in M.
\end{equation}
A Riemannian metric on $M$ is said to be {\it adapted} to the
dominated splitting if we can take $C=1$ in \eqref{domi} with
respect to this metric. Adapted metric always exists, see \cite{Go}
for details.

Although restriction of the dominated splitting is much weaker than
(partially) hyperbolic splitting, we show the restriction is strong
enough to exclude the possibility of minimality. Recall that the map
$f$ is said to be minimal if for each $x\in M$, the orbit
$\mathcal{O}_f(x)=\{f^nx:n\in\mathbb{Z}\}$ is a dense subset in $M$.
The following is our main result.

\vskip.4cm

\noindent{\bf Main Theorem.} {\it Let $M$ be a closed Riemannian
manifold and $f:M\to M$ be a diffeomorphism on $M$. If $f$ has a
global dominated splitting, then $f$ can not be minimal.}

\vskip.4cm

There are vast results in the case that the dimension of $M$ is $2$.
Pujals and Sambarino gave several good characteristics in
\cite{PS1,PS2} for the topology of invariant set having some
dominated splitting for $C^2$ diffeomorphisms on surfaces. On the
other hand, Xia proved in \cite{X} that for all compact surface $M$
with nonzero Euler characteristic, any homeomorphism on $M$ admits
some periodic points. In \cite{Z} X. Zhang considered a
diffeomorphism $f$ on a closed surface and an $f$-invariant set
$\Lambda$ with dominated splitting. He used Liao's selecting lemma
and Crovisier's central models to find a periodic orbit near
$\Lambda$. Our result follows from an observation from \cite{Z}
combining with Ma\~n\'e's argument and Liao's sifting lemma (or use
Liao's selecting lemma).

\section{partial hyperbolicity and quasi-hyperbolic strings}

In this section let's review several useful results for later
discussions. Let $M$ be a compact manifold and assume $f:M\to M$ has
a dominated splitting $TM=E\oplus F$. We always assume the
Riemannian metric on $M$ is chosen to be {\it adapted}. That is,
there exists $\lambda\in(0,1)$ with
$$\|Df|_{E(x)}\|\cdot\|Df^{-1}|_{F(fx)}\|<\lambda,\text{ for all }x\in M.$$
In this case we say $TM=E\oplus F$ is a $\lambda$-dominated
splitting. Generally we have
$\|Df^n|_{E(x)}\|\le\prod_{k=0}^{n-1}\|Df|_{E(f^kx)}\|$ for all
$n\ge1$. We have similar observations for the subbundle $F$.

The first result we recall is the argument of Ma\~n\'e which can
locate a {\it nonrecurrent} point (see \cite[Lemma 5.2]{M}. Also see
\cite[Corollary 1]{BRHRHTU}).
\begin{pro}\label{mane}
Let $f$ be a diffeomorphism on $M$ and $\mathcal{W}$ be an
$f$-invariant foliation tangent to a distribution $E\subset TM$ such
that $Df$ is uniformly expanding (or uniformly contracting) on $E$.
Then there exists a nonrecurrent point of $f$. Moreover the set
$\{z\in M: z\notin\omega(z)\}$ of points that are nonrecurrent in
the future is dense in every leaf of $\mathcal{W}$.
\end{pro}
Let's sketch the proof here. First we show that each leaf $W(x)$
contains at most one periodic point. To this end, we assume $E$ is
$\nu$-expanding for some $\nu>1$. Suppose that there are two
periodic points $p,q\in W(x)$ for some $x\in M$. Pick some $n\ge1$
with $f^np=p$ with $f^nq=q$. Let $\gamma$ be a smooth curve in
$W(x)$ connecting $p$ and $q$ with length $|\gamma|\le
\nu^{1/2}\cdot d_{W}(p,q)$. Then $f^{-n}\gamma$ is also a path
connecting $p$ and $q$ with length $d_W(p,q)\le |f^{-n}\gamma|\le
\nu^{-n}\cdot|\gamma|\le\nu^{-1/2}\cdot d_W(p,q)$. This is
impossible unless $p=q$. This shows that each leaf $W(x)$ contains
at most one periodic point. Let $N\ge1$ such that $\lambda^N\ge5$.
Then for a nonperiodic point $y\in W(x)$ we pick $\epsilon>0$ small
enough such that $f^kB(y,\epsilon)$, $0\le k\le N$ are pairwisely
disjoint. Then choose $\delta>0$ to be much smaller than $\epsilon$.
We define inductively a sequence of closed disks $D_n\subset
f^nW(y,\delta)$ for $n\ge1$ with $D_{n+1}\subset fD_n$ and $D_n\cap
B(y,2\delta)=\emptyset$. Then each point $z\in
\bigcap_{n\ge1}f^{-n}D_n$ satisfies $f^nz\notin B(y,2\delta)$ for
all $n\ge1$. Such $z$ is nonrecurrent and we finishes the proof. See
\cite[Lemma 5.2]{BRHRHTU} for the construction of $D_n$ and more
details.

The second one is Liao's {\it Sifting Lemma} which helps us to
locate some {\it periodic} point. See \cite{L1,L2} for example.
\begin{pro}\label{liao}
Let $1\le I\le d-1$ and $\Lambda$ be a compact invariant set of $f$
with a $\lambda$-dominated splitting $T_{\Lambda}M=E\oplus F$ of
index $I$. Assume
\begin{enumerate}
\item There is a point $b\in\Lambda$ satisfying
$\prod_{k=0}^{n-1}\|Df|_{E(f^kb)}\|\ge1$ for all $n\ge1$.

\item (The tilde condition.) There are $\lambda_1$ and $\lambda_2$ with
$\lambda<\lambda_1<\lambda_2<1$ such that if a point $x\in\Lambda$
satisfies $\prod_{k=0}^{n-1}\|Df|_{E(f^kx)}\|\ge\lambda_2^n$ for all
$n\ge1$, then the omega-set $\omega(x)$ contains a point $c$
satisfying $\prod_{k=0}^{n-1}\|Df|_{E(f^kc)}\|\le\lambda_1^n$ for
all $n\ge1$.
\end{enumerate}

Then for each $\lambda_3\in(\lambda_2,1)$ and each $l\in\mathbb{N}$,
there are $l$ positive integers $n_1<n_2<\cdots<n_l$ with the
following property: for every $j=1,\cdots,l-1$ and every
$k=n_j+1,\cdots,n_{j+1}$,
\begin{equation}\label{double}
\prod_{i=n_j}^{k-1}\|Df|_{E(f^ib)}\|\le\lambda_3^{k-n_j}\text{ and }
\prod_{i=k-1}^{n_{j+1}-1}\|Df|_{E(f^ib)}\|\ge\lambda_2^{n_{j+1}-k+1}.
\end{equation}
\end{pro}
See \cite[Lemma 2.2]{W} for a proof. In the following we sketch how
to find a hyperbolic periodic point near $\Lambda$. Proposition
\ref{liao} shows that there are many `double' uniform stings when
$f^nb$ approaches some `good' point $c\in\omega(b)$: for each
$\lambda_3\in(\lambda_2,1)$ and each $l\in\mathbb{N}$, there are $l$
positive integers $n_1<n_2<\cdots<n_l$ with the following property:
for every $j=1,\cdots,l-1$ and every $k=n_j+1,\cdots,n_{j+1}$,
\begin{equation*}
\prod_{i=n_j}^{k-1}\|Df|_{E(f^ib)}\|\le\lambda_3^{k-n_j}\text{ and }
\prod_{i=k-1}^{n_{j+1}-1}\|Df|_{E(f^ib)}\|\ge\lambda_2^{n_{j+1}-k+1}.
\end{equation*}
Then by $\lambda$-domination assumption we have that
$$\prod_{i=k-1}^{n_{j+1}-1}\|Df^{-1}|_{F(f^{i+1}b)}\|\le
\prod_{i=k-1}^{n_{j+1}-1}\frac{\lambda}{\|Df|_{E(f^ib)}\|}
\le\left(\frac{\lambda}{\lambda_2}\right)^{n_{j+1}-k+1}$$ for every
$j=1,\cdots,l-1$ and every $k=n_j+1,\cdots,n_{j+1}$. Let
$\tilde{\lambda}=\max\{\sqrt{\lambda},\lambda_3,\lambda/\lambda_2\}$.
Then $\tilde{\lambda}<1$ and $(f^{n_j}b,f^{n_{j+1}}b)$ forms a
`$\tilde{\lambda}$-quasi-hyperbolic string' for each
$j=1,\cdots,l-1$.

Let $L\ge1$ and $d_0$ given by \cite[Theorem 1.1]{G} with respect to
$\tilde{\lambda}$. For $\epsilon\in(0,d_0]$ let's pick an integer
$l=l(\epsilon)\ge1$ large enough such that given arbitrary $l$
points $x_1,\cdots,x_l$ in $M$, there exists $1\le i<j\le l$ such
that $d(x_i,x_j)<\epsilon$. For this $l$ we let $n_1<\cdots<n_l$ be
given by Proposition \ref{liao} such that \eqref{double} holds. Then
$d(f^{n_i}b,f^{n_j}b)<\epsilon$ for some $1\le i<j\le l$. Finally we
apply Liao--Gan's shadowing lemma (see \cite[Theorem 1.1]{G}) to
find a hyperbolic periodic point $L\epsilon$-shadowing the periodic
pseudo-orbit
$\{(f^{n_i}b,f^{n_{i+1}}b),(f^{n_{i+1}}b,f^{n_{i+2}}b),\cdots,(f^{n_{j-1}}b,f^{n_j}b)\}$.
This finishes the proof. For more details see \cite{G,W}. Also see
Liao's Selecting Lemma \cite[Lemma 2.3]{W} for more information.

\section{dominated splitting and minimality}

With the preparations in previous section let's prove the main
theorem that if the map $f$ has a global dominated splitting, then
it can not be minimal.

\vskip.4cm

\noindent{\it Proof of Main Theorem.} Let $TM=E\oplus F$ be a
$\lambda$-dominated splitting for some $\lambda\in(0,1)$ with
respect some adapted Riemannian metric. If $F$ is uniformly
expanding, then $F$ is uniquely integrable and tangent to the
strongly unstable foliation $\mathcal{W}^{su}$. By Proposition
\ref{mane} there exists some nonrecurrent point of $f$ and hence $f$
can not be minimal on $M$. Similarly if $E$ is uniformly
contracting, then $f$ can not be minimal either. Then we are left
with the case that neither $E$ is uniformly contracting, nor $F$ is
uniformly expanding. In this case we have that:
\begin{enumerate}
\item since $E$ is not uniformly contracting, there exists
$p\in M$ such that $\|Df^n|_{E(p)}\|\ge1$ for all $n\ge1$,

\item since $F$ is not uniformly expanding, there exists $q\in M$ such that
$\|Df^{-n}|_{F(f^nq)}\|\ge1$ for all $n\ge1$.
\end{enumerate}
We first observe that, by the $\lambda$-domination assumption, for
all $n\ge1$,
\begin{equation}\label{good}
\prod_{k=0}^{n-1}\|Df|_{E(f^kq)}\|
\le\prod_{k=0}^{n-1}\frac{\lambda}{\|Df^{-1}|_{F(f^{k+1}q)}\|}
\le\frac{\lambda^n}{\|Df^{-n}|_{F(f^nq)}\|}\le\lambda^n.
\end{equation}
Also note that the first condition in Proposition \ref{liao} is
already satisfied if we take $\Lambda=M$ and $b=p$ since
$\prod_{k=0}^{n-1}\|Df|_{E((f^kp)}\|\ge\|Df^n|_{E(p)}\|\ge1$ for
each $n\ge1$. Then we divide the discussion into two subcases:

Subcase 1. The {\it tilde condition} in Proposition \ref{liao} holds
on $M$ for some $\lambda_1,\lambda_2$ with
$\lambda<\lambda_1<\lambda_2<1$. Then by Proposition \ref{liao} and
Liao--Gan's shadowing lemma (Theorem 1.1 in \cite{G}), there does
exist a hyperbolic periodic point of $f$: the map $f$ can not be
minimal.

Subcase 2. The {\it tilde condition} fails. So for each pair
$\lambda_1,\lambda_2$ with $\lambda<\lambda_1<\lambda_2<1$, there
exists some point $\tilde{x}\in M$ such that
\begin{itemize}
\item $\prod_{k=0}^{n-1}\|Df|_{E(f^k\tilde{x})}\|\ge\lambda_2^n$ for all
$n\ge1$.

\item for each $y\in\omega(\tilde{x})$, there exists some $n(y)\ge1$ with
$\prod_{k=0}^{n(y)-1}\|Df|_{E(f^ky)}\|\ge\lambda_1^{n(y)}$.
\end{itemize}
According to \eqref{good}, we see that $q\notin \omega(\tilde{x})$
since $\lambda<\lambda_1$. So $\omega(\tilde{x})\subsetneq M$ for
some point $\tilde{x}\in M$ and the map $f$ is not minimal either.

This finishes the verification for both subcases and ends the proof
of theorem.

\begin{remark}
The result is not true if we consider invariant subsets instead of
the whole manifold, since there are various kinds of minimal subsets
on which the map $f$ is dominated. For example let
$A=\begin{pmatrix} 2 & 1
\\ 1 & 1
\end{pmatrix}$ and $f_A:\mathbb{T}^2\to\mathbb{T}^2$ be the induced
diffeomorphism. Let $R:\mathbb{T}\to\mathbb{T}$ be an irrational
rotation. Then $f_A$ is Anosov with a fixed point $o\in\mathbb{T}^2$
and $R$ is minimal. Moreover the product system
$(R,f_A):\mathbb{T}\times\mathbb{T}^2\to\mathbb{T}\times\mathbb{T}^2$
is partially hyperbolic with an invariant minimal subset
$\Lambda=\mathbb{T}\times\{o\}$.
\end{remark}

\section*{Acknowledgments} This work is done when the author was visiting
Peking University. We would like to thank Shaobo Gan, Lan Wen and
Zhihong Xia for discussions and suggestions. Especially we thank
Baolin He and Peng Sun for many helpful explanations.

\end{document}